# CONSISTENT MASS MATRIX OF TEN NODES TETRAHEDRAL ELEMENT BASED ON ANALYTICAL INTEGRATION.

E Hanukah

Faculty of Mechanical Engineering, Technion – Israel Institute of Technology, Haifa 32000, Israel
Corresponding author Email: eliezerh@tx.technion.ac.il

**Abstract**

Currently, components of consistent mass matrix are computed using various numerical integration schemes, each one alters in number of integration (Gauss) points, requires different amount of computations and possess different level of accuracy. We discuss the closed-form mass matrix based on analytical integration. Curved-sided and straight-sided elements are considered. For a straight-sided element we derive an exact analytical easy to implement consistent mass matrix. For a curved-sided element an exact analytical mass matrix is derived, however it is rather lengthy, hence approximations are proposed. Three systematic approximations to the metric (jacobian determinant) are suggested; constant metric (CM), linearly varying metric (LM) and quadratic metric (QM). CM requires evaluation of the metric at the centroid, LM requires metric evaluations at the four corner nodes and QM uses metric values at all the ten nodes. Analytical integration together with approximated metric models yields closed-form semi-analytical mass matrices. The accuracy of the schemes is studied numerically using randomly generated coarse mesh. Our findings reveal significant superiority in accuracy and computations over equivalent schemes. An important implication of this study is that based on the results, it is superior to use our CM, LM and QM semi-analytical mass matrices over mass matrices based on numerical integration schemes which involve four, five and fifteen point Gauss quadrature. For a straight-sided element, CM, LM and QM admit an exact consistent mass matrix.

## 1. Introduction

Automatic mesh generators are widely used in finite element analysis. Complicated three dimensional geometries are often meshed with tetrahedral elements, in fact, it is easier to produce tetrahedral than hexahedral mesh (e.g. [1-3]). Four node tetrahedral elements admit only homogeneous deformation and well known for their stiff behavior and volumetric locking, therefore they are not recommended. Hence, ten node tetrahedral elements are of great importance.

Mass matrix components, internal forces and stiffness matrix, all require integration in the element domain, which is most commonly obtained with the help on numerical integration schemes e.g. [4, 5]. Several studies exist that exploit the idea of closed-form integration for stiffness matrixes [6-10], significant time savings is established. Furthermore, hierarchical semi-analytical displacement based approach is used to model three dimensional finite bodies e.g. [11-14] yielding new analytical solutions.

In present study we follow the basic guidelines presented in [15]. Analytical integration together with symbolic approximations is used for computation of consistent mass matrix. Straight-sided (straight edges, flat faces) element is characterized by constant metric namely jacobian determinant of global-local coordinate system is independent of natural coordinates. Analytic integration result in simple, easy to implement expressions.

For a curved-sided element, jacobian matrix is linear and the metric is cubic with respect to coordinates. With the help of Taylor's multivariable expansion, convenient representation of the metric is derived. Analytical integration result in exact mass matrix, which is used later in numerical study as reference values for error computation for approximate schemes. Explicit closed-form representation of an exact mass matrix components is rather lengthy and not easy to implement, therefore for practical applications three systematic approximations are suggested; constant metric (CM), linear metric (LM) and quadratic metric (QM). CM equal to the metric evaluated at the centroid of the element, LM requires evaluation of the metric at corner nodes and linear ansatz functions while QM involves evaluation at all the nodes and element shape functions. Analytic integration is then applied and approximated mass matrices are derived.

CM, LM and QM are exact for arbitrary configuration of straight sided element. Preliminary numerical study is conducted to test performance of new CM, LM and QM mass matrices. Randomly generated coarse mesh is produced and an averaged absolute error is calculated with respect to exact results. In terms of computations, CM is equivalent to one point scheme; however it significantly over-performs even four and five point numerical integration schemes. QM and fifteen point numerical integration scheme demonstrated similar accuracy; nevertheless, QM requires ten metric evaluations while fifteen evaluations are needed for numerical integrations.

The outline of the paper is as follows. Section 2 recalls all the details of consistent mass matrix formulation, natural coordinates, shape functions, metric and jacobian matrix definitions, numerical and analytical integration in the element domain. Section 3 record details of derivation of an exact analytical mass matrix for straight-sided (constant metric) element. It is shown than the resulting exact mass matrix computationally equivalent to numerical integration using one integration point. Section 4 considers a curved sided (varying metric) element. Representation of the metric is developed using Taylor's multivariable expansion about the origin. An exact analytical mass matrix is derived. Three systematic approximations for the metric are suggested: constant metric (CM), linear metric (LM) and quadratic metric (QM). An analytical integration is used to derive closed-form semi-analytical mass matrices. Section 5 contains preliminary numerical accuracy study and its comparison to equivalent numerical integration schemes. It is found that our CM mass matrix over-performs in accuracy and efficiency numerically integrated mass matrix using four and five Gauss points. In addition it is found that numerically integrated mass matrix based on 15 integration points requires about 50% more computation than our QM mass matrix, though the same accuracy is established. Section 6 records our conclusions.

## 2. Background

Initial nodal locations of the standard 10 node tetrahedral element (e.g. [16] pp.72, [17] pp.120) are denoted by $\mathbf{X}_i \, (i=1,..,10)$, where components are given in terms of global Cartesian coordinates system $\mathbf{X}_i = X_{ki}\mathbf{e}_k \, (k=1,2,3, i=1,..,10)$, summation convention on repeated index is implied. Here and throughout the text, bold symbols traditionally denote vector or tensor quantities. Local convected coordinate system (natural coordinates) $\{\xi,\eta,\zeta\}$ admits

$$0 \le \xi \le 1-\eta-\zeta \, , \, 0 \le \eta \le 1-\zeta \, , \, 0 \le \zeta \le 1 \tag{1}$$

The standard shape functions $\varphi^i \, (i=1,..,10)$ in terms of natural coordinates is given by

$$\begin{aligned}
\varphi^1 &= (1-\xi-\eta-\zeta)(1-2\xi-2\eta-2\zeta) \quad , \quad \varphi^2 = \xi(2\xi-1) \\
\varphi^3 &= \eta(2\eta-1) \, , \, \varphi^4 = \zeta(2\zeta-1) \, , \, \varphi^5 = 4\xi(1-\xi-\eta-\zeta) \\
\varphi^6 &= 4\xi\eta \, , \, \varphi^7 = 4\eta(1-\xi-\eta-\zeta) \, , \, \varphi^8 = 4\zeta(1-\xi-\eta-\zeta) \\
\varphi^9 &= 4\xi\zeta \quad , \quad \varphi^{10} = 4\eta\zeta
\end{aligned} \tag{2}$$

Material point X occupies location $\mathbf{X}$ inside the element domain (1) is given by

$$\mathbf{X} = \varphi^i \mathbf{X}_i \quad (i=1,..,10) \tag{3}$$

In this study we consider homogeneous initial configuration $\rho_0 = \text{const}$. Extension for linearly varying initial density $\rho_0 = (1-\xi-\eta-\zeta)\rho_1 + \xi\rho_2 + \eta\rho_3 + \zeta\rho_4$ or quadratic initial density $\rho_0 = \varphi^i \rho_i \, (i=1,..,10)$ etc., where $\rho_i$ denote initial nodal densities, follows the same steps.
The metric or jacobian determinant of global-local coordinates transformation - J is given by

$$J = \mathbf{X}_{,1} \times \mathbf{X}_{,2} \cdot \mathbf{X}_{,3} = \begin{vmatrix} (\mathbf{X}\cdot\mathbf{e}_1)_{,1} & (\mathbf{X}\cdot\mathbf{e}_1)_{,2} & (\mathbf{X}\cdot\mathbf{e}_1)_{,3} \\ (\mathbf{X}\cdot\mathbf{e}_2)_{,1} & (\mathbf{X}\cdot\mathbf{e}_2)_{,2} & (\mathbf{X}\cdot\mathbf{e}_2)_{,3} \\ (\mathbf{X}\cdot\mathbf{e}_3)_{,1} & (\mathbf{X}\cdot\mathbf{e}_3)_{,2} & (\mathbf{X}\cdot\mathbf{e}_3)_{,3} \end{vmatrix} > 0 \tag{4}$$

$$J_{mn}(\xi,\eta,\zeta,X_{ki}) = (\mathbf{X}\cdot\mathbf{e}_m)_{,n} \, , \, (i=1,..,10, m,n,k=1,2,3)$$

Where $(\times)$ and $(\cdot)$ stand for vector cross and scalar products and $|\cdot|$ stand for determinant operator, comma denotes partial differentiation with respect to coordinates. Here and throughout the text, determinant of general (non-symmetric) 3x3 matrixes is computed using the next consistent with standard determinant definition formula

$$J = J_{11}J_{22}J_{33} - J_{11}J_{23}J_{32} - J_{31}J_{22}J_{13} - J_{21}J_{12}J_{33} + J_{21}J_{32}J_{13} + J_{31}J_{12}J_{23} \quad (5)$$

Differential volume element $dV$, and initial volume $V$ are defined by (e.g.[18])

$$dV = J d\xi d\eta d\zeta \ , \ V = \int_V dV = \int_0^{+1} \int_0^{1-\zeta} \int_0^{1-\eta-\zeta} J d\xi d\eta d\zeta \quad (6)$$

Isoparametric formulation (e.g.[17] pp.104) for mass conserving element, lead to the next consistent, symmetric, and positive definite mass matrix

$$M^{ij} = \int_V \rho_0 \phi^i \phi^j dV \ , \ M^{ij} = M^{ji} \ , \ (i,j = 1,..,10) \quad (7)$$

Standard numerical integration in an element domain is recalled (e.g. [17] pp.122, [16] pp.79)

$$\int_V f(\xi,\eta,\zeta) dV = \sum_{p=1}^{n_p} f(\xi_p,\eta_p,\zeta_p) J(\xi_p,\eta_p,\zeta_p) w_p \quad (8)$$

Where $n_p$ stand for number of integration (Gauss) points, $w_p$ denotes weights and $\xi_p, \eta_p, \zeta_p$ are coordinates of integration points. Widely used five point quadrature is detailed at [17] pp.122, four and fifteen point schemes are detailed in [16] pp.79. For later convenience, numerical computation of the mass matrix (7)(8) represented in more details

$$n_p = 1 \ , \ M^{ij} = \rho_0 \hat{J}_1 w_1 \hat{M}_1^{ij}$$

$$n_p = 4 \ , \ M^{ij} = \rho_0 (\hat{J}_1 w_1 \hat{M}_1^{ij} + \hat{J}_2 w_2 \hat{M}_2^{ij} + \hat{J}_3 w_3 \hat{M}_3^{ij} + \hat{J}_4 w_4 \hat{M}_4^{ij})$$

$$n_p = 5 \ , \ M^{ij} = \rho_0 (\hat{J}_1 w_1 \hat{M}_1^{ij} + \hat{J}_2 w_2 \hat{M}_2^{ij} + \hat{J}_3 w_3 \hat{M}_3^{ij} + \hat{J}_4 w_4 \hat{M}_4^{ij} + \hat{J}_5 w_5 \hat{M}_5^{ij}) \quad (9)$$

$$n_p = 15 \ , \ M^{ij} = \rho_0 (\hat{J}_1 w_1 \hat{M}_1^{ij} + ... + \hat{J}_{10} w_{10} \hat{M}_{10}^{ij} + ... + \rho_0 \hat{J}_{15} w_{15} \hat{M}_{15}^{ij})$$

$$\hat{M}_p^{ij} = \varphi_p^i \varphi_p^j \ (p = 1,..,n_p, i,j = 1,..,10)$$

Where $\hat{J}_p \ (p = 1,..,n_p)$ stand for the metric evaluation at the integration point $p$, $\varphi_p^i$ denote shape function $i$ evaluated at integration point $p$.

## 3. Straight-sided element.

It is always possible to generate initial mesh using straight-sided (straight-edges / flat faces / constant metric) tetrahedral elements; nodes 5-10 are located in the middle of their edges

$$\mathbf{X}_5 = \frac{1}{2}(\mathbf{X}_1 + \mathbf{X}_2) \ , \ \mathbf{X}_6 = \frac{1}{2}(\mathbf{X}_2 + \mathbf{X}_3) \ , \ \mathbf{X}_7 = \frac{1}{2}(\mathbf{X}_1 + \mathbf{X}_3)$$
$$\mathbf{X}_8 = \frac{1}{2}(\mathbf{X}_1 + \mathbf{X}_4) \ , \ \mathbf{X}_9 = \frac{1}{2}(\mathbf{X}_2 + \mathbf{X}_4) \ , \ \mathbf{X}_{10} = \frac{1}{2}(\mathbf{X}_3 + \mathbf{X}_4) \tag{10}$$

Using (3),(4) and the above relations (10), components of the jacobian matrix become independent of the coordinates and are given by

$$J_{mn} = \begin{pmatrix} X_{12} - X_{11} & -X_{11} + X_{13} & -X_{11} + X_{14} \\ X_{22} - X_{21} & -X_{21} + X_{23} & -X_{21} + X_{24} \\ X_{32} - X_{31} & -X_{31} + X_{33} & -X_{31} + X_{34} \end{pmatrix} \tag{11}$$

Using the determinant formula (5) and the above jacobian matrix it follows that the metric is constant. Analytical integration (7) return exact consistent mass matrix

$$M^{ij} = \rho_0 \frac{J}{2520} M_0^{ij} \ , \ (i,j = 1,..,10) \tag{12}$$

Where $M_0^{ij}$ is given by

$$M_0^{ij} = \begin{pmatrix} 6 & 1 & 1 & 1 & -4 & -6 & -4 & -4 & -6 & -6 \\ 1 & 6 & 1 & 1 & -4 & -4 & -6 & -6 & -4 & -6 \\ 1 & 1 & 6 & 1 & -6 & -4 & -4 & -6 & -6 & -4 \\ 1 & 1 & 1 & 6 & -6 & -6 & -6 & -4 & -4 & -4 \\ -4 & -4 & -6 & -6 & 32 & 16 & 16 & 16 & 16 & 8 \\ -6 & -4 & -4 & -6 & 16 & 32 & 16 & 8 & 16 & 16 \\ -4 & -6 & -4 & -6 & 16 & 16 & 32 & 16 & 8 & 16 \\ -4 & -6 & -6 & -4 & 16 & 8 & 16 & 32 & 16 & 16 \\ -6 & -4 & -6 & -4 & 16 & 16 & 8 & 16 & 32 & 16 \\ -6 & -6 & -4 & -4 & 8 & 16 & 16 & 16 & 16 & 32 \end{pmatrix} \tag{13}$$

It is important to emphasize, that (12) has the same form as numerical integration using one Gauss point (9), as a result, computational effort associated with (12) is similar. On the other hand, while our mass matrix (12) is exact, mass matrix built on one point numerical integration is unacceptably poor approximation, it will be shown later that even four and five point numerical integration is inaccurate for straight sided elements.

Here and throughout the study, computer algebra system (CAS) MAPLE[TM] were used to perform all the symbolic manipulations, including integration, differentiation, simplification, direct translation of explicit expression to Fortran77, generation of pseudo-random numbers for coarse mesh etc.

## 4. Curved-sided element.

For a general configuration of ten nodes tetrahedral element, namely nodes 5-10 are not necessarily in the middle of the edges, jacobian matrix (4) is represented by

$$J_{mn} = J^0_{mn} + \xi J^1_{mn} + \eta J^2_{mn} + \zeta J^3_{mn} \ , \ (m,n = 1,2,3) \tag{14}$$

Where $J^k_{mn}$ (k = 0,..,3) are functions of nodal components $X_{ki}$ (k = 1,2,3, i = 1,..,10) only

$$J^0_{mn} = \begin{pmatrix} -X_{12} - 3X_{11} + 4X_{15} & 4X_{17} - X_{13} - 3X_{11} & -X_{14} - 3X_{11} + 4X_{18} \\ -X_{22} - 3X_{21} + 4X_{25} & 4X_{27} - X_{23} - 3X_{21} & -X_{24} - 3X_{21} + 4X_{28} \\ -X_{32} - 3X_{31} + 4X_{35} & 4X_{37} - X_{33} - 3X_{31} & -X_{34} - 3X_{31} + 4X_{38} \end{pmatrix}$$

$$J^1_{mn} = 4\begin{pmatrix} X_{11} + X_{12} - 2X_{15} & -X_{17} + X_{16} - X_{15} + X_{11} & -X_{15} + X_{11} - X_{18} + X_{19} \\ X_{21} + X_{22} - 2X_{25} & -X_{27} + X_{26} - X_{25} + X_{21} & -X_{25} + X_{21} - X_{28} + X_{29} \\ X_{31} + X_{32} - 2X_{35} & -X_{37} + X_{36} - X_{35} + X_{31} & -X_{35} + X_{31} - X_{38} + X_{39} \end{pmatrix}$$

$$J^2_{mn} = 4\begin{pmatrix} -X_{17} + X_{16} - X_{15} + X_{11} & -2X_{17} + X_{13} + X_{11} & X_{110} + X_{11} - X_{18} - X_{17} \\ -X_{27} + X_{26} - X_{25} + X_{21} & -2X_{27} + X_{23} + X_{21} & X_{210} + X_{21} - X_{28} - X_{27} \\ -X_{37} + X_{36} - X_{35} + X_{31} & -2X_{37} + X_{33} + X_{31} & X_{310} + X_{31} - X_{38} - X_{37} \end{pmatrix} \tag{15}$$

$$J^3_{mn} = 4\begin{pmatrix} -X_{15} + X_{11} - X_{18} + X_{19} & X_{110} + X_{11} - X_{18} - X_{17} & X_{14} + X_{11} - 2X_{18} \\ -X_{25} + X_{21} - X_{28} + X_{29} & X_{210} + X_{21} - X_{28} - X_{27} & X_{24} + X_{21} - 2X_{28} \\ -X_{35} + X_{31} - X_{38} + X_{39} & X_{310} + X_{31} - X_{38} - X_{37} & X_{34} + X_{31} - 2X_{38} \end{pmatrix}$$

With the help of (4),(5) and Taylor's multivariable expansion about the origin $\mathbf{X}_1 = \mathbf{X}(0,0,0)$, the metric J is represented as

$$\begin{aligned} J = J_0 &+ \\ & \xi J_1 + \eta J_2 + \zeta J_3 + \\ & \xi\eta J_4 + \xi\zeta J_5 + \eta\zeta J_6 + \xi^2 J_7 + \eta^2 J_8 + \zeta^3 J_9 + \\ & \xi\xi\eta J_{10} + \xi\eta\eta J_{11} + \xi\eta\zeta J_{12} + \xi^2\zeta J_{13} + \eta^2\zeta J_{14} + \xi\zeta^2 J_{15} + \eta\zeta^2 J_{16} + \xi^3 J_{17} + \eta^3 J_{18} + \zeta^3 J_{19} \end{aligned} \tag{16}$$

Where $J_w$ (w = 0,..,19) depend on nodal components $X_{ki}$ (k = 1,2,3, i = 1,..,10) but independent of natural coordinates $\{\xi, \eta, \zeta\}$

$$\begin{aligned} J_0 &= J(\mathbf{X}_1) \\ J_1 &= \frac{\partial J(\mathbf{X}_1)}{\partial \xi} \ , \ J_2 = \frac{\partial J(\mathbf{X}_1)}{\partial \eta} \ , \ J_3 = \frac{\partial J(\mathbf{X}_1)}{\partial \zeta} \\ J_4 &= \frac{\partial^2 J(\mathbf{X}_1)}{\partial \xi \partial \eta} \ , \ .. \ , \ J_{12} = \frac{\partial^3 J(\mathbf{X}_1)}{\partial \zeta \partial \eta \partial \zeta} \ , \ .. \ , \ J_{19} = \frac{\partial^3 J(\mathbf{X}_1)}{\partial \zeta^3} \\ J_w &= \breve{J}_w(X_{ki}) \ , \ (w = 0,..,19, k = 1,2,3, i = 1,..,10) \end{aligned} \tag{17}$$

Using the above representation together with mass matrix definition (7), shape function definition (2) and volume integration (6), exact mass matrix components are computed and used as an exact values in later numerical study. One of the exact terms is given to illustrate the general form

$$M^{55}=\frac{\rho_0}{56700}(720J_0+270J_1+90J_2+90J_3+120J_4+30J_5+20J_6+30J_7+10J_8+20J_9+60J_{10}+ \\ 12J_{11}+6J_{12}+6J_{13}+12J_{14}+3J_{15}+2J_{16}+6J_{17}+2J_{18}+6J_{19}) \quad (18)$$

Explicit expressions for $J_w$ $(w=0,..,19)$ as a function of $X_{ki}$ $(k=1,2,3, i=1,..,10)$ and explicit expressions to all the exact analytical mass matrix components $M^{ij}$ $(i,j=1,..,10)$ as a function of $J_w$ $(w=0,..,19)$ are omitted from the text for a sake of brevity. Please contact the corresponding author if needed. Though an exact computation of the consistent mass matrix is performed, it is rather computationally expensive since $J_w(X_{kI})$ $(w=0,..,19)$ are relatively lengthy expressions, therefore simplified models for the metric $J$ are suggested.

The simplest approximation for the metric is a constant metric - CM.

$$J \approx J_{cent} \ , \ J_{cent} = \left|J_{mn}^{cent}\right| \quad (19)$$

Where $J_{mn}^{cent}$ stand for jacobian matrix (14) evaluated at the centroid $\mathbf{X}_{cent} = \mathbf{X}(\frac{1}{4},\frac{1}{4},\frac{1}{4})$

$$J_{mn}^{cent} = \begin{pmatrix} -X_{17}+X_{16}-X_{18}+X_{19} & X_{16}-X_{15}+X_{110}-X_{18} & -X_{15}+X_{19}+X_{110}-X_{17} \\ -X_{27}+X_{26}-X_{28}+X_{29} & X_{26}-X_{25}+X_{210}-X_{28} & -X_{25}+X_{29}+X_{210}-X_{27} \\ -X_{37}+X_{36}-X_{38}+X_{39} & X_{36}-X_{35}-X_{38}+X_{310} & -X_{35}+X_{39}-X_{37}+X_{310} \end{pmatrix} \quad (20)$$

CM approximation (19) together with analytical integration (6) and definition of $M_0^{ij}$ given by (13) yield CM mass matrix for curved-sided element

$$M_{CM}^{ij} = \rho_0 \frac{J_0}{2520} M_0^{ij} \ , \ (i,j=1,..,10) \quad (21)$$

The above is in agreement with the exact mass matrix of a straight-sided element (12), namely for a straight-sided element the above (21) equal to the exact (12).

Our second model for the metric is linearly varying metric - LM namely

$$J \approx (1-\xi-\eta-\zeta)\tilde{J}_1 + \xi\tilde{J}_2 + \eta\tilde{J}_3 + \zeta\tilde{J}_4 \\ \tilde{J}_k = \left|\tilde{J}_{mn}^k\right| \ , \ \tilde{J}_{mn}^k = J_{mn}\big|_{node\,k} \ , \ (k=1,2,3,4) \quad (22)$$

Where $\tilde{J}_k$ stand for metric evaluated at the node $k=1,2,3,4$, and $\tilde{J}_{mn}^k$ is the jacobian matrix (14) at the node $k=1,2,3,4$.

$$\tilde{J}^1_{mn} = \begin{pmatrix} -X_{12}-3X_{11}+4X_{15} & 4X_{17}-X_{13}-3X_{11} & -X_{14}-3X_{11}+4X_{18} \\ -X_{22}-3X_{21}+4X_{25} & 4X_{27}-X_{23}-3X_{21} & -X_{24}-3X_{21}+4X_{28} \\ -X_{32}-3X_{31}+4X_{35} & 4X_{37}-X_{33}-3X_{31} & -X_{34}-3X_{31}+4X_{38} \end{pmatrix}$$

$$\tilde{J}^2_{mn} = \begin{pmatrix} X_{11}+3X_{12}-4X_{15} & 4X_{16}-4X_{13}-X_{11}-4X_{13} & -4X_{15}-X_{11}+4X_{19}-X_{14} \\ X_{21}+3X_{22}-4X_{25} & 4X_{26}-4X_{23}-X_{21}-4X_{23} & -4X_{25}-X_{21}+4X_{29}-X_{24} \\ X_{31}+3X_{32}-4X_{35} & 4X_{36}-4X_{33}-X_{31}-4X_{33} & -4X_{35}-X_{31}+4X_{39}-X_{34} \end{pmatrix}$$

$$\tilde{J}^3_{mn} = \begin{pmatrix} -4X_{17}+4X_{16}+X_{11}-X_{12} & X_{11}+3X_{13}-4X_{17} & 4X_{110}+X_{11}-4X_{17}-X_{14} \\ -4X_{27}+4X_{26}+X_{21}-X_{22} & X_{21}+3X_{23}-4X_{27} & 4X_{210}+X_{21}-4X_{27}-X_{24} \\ -4X_{37}+4X_{36}+X_{31}-X_{32} & X_{31}+3X_{33}-4X_{37} & 4X_{310}+X_{31}-4X_{37}-X_{34} \end{pmatrix}$$

$$\tilde{J}^4_{mn} = \begin{pmatrix} X_{11}-4X_{18}+4X_{19}-X_{12} & 4X_{110}+X_{11}-4X_{18}-X_{13} & 3X_{14}+X_{11}-4X_{18} \\ X_{21}-4X_{28}+4X_{29}-X_{22} & 4X_{210}+X_{21}-4X_{28}-X_{23} & 3X_{24}+X_{21}-4X_{28} \\ X_{31}-4X_{38}+4X_{39}-X_{32} & 4X_{310}+X_{31}-4X_{38}-X_{33} & 3X_{34}+X_{31}-4X_{38} \end{pmatrix}$$

(23)

LM approximation (22) together with shape functions definition (2) and analytical integration (6) applied to mass matrix definition (7) result in

$$M^{ij}_{LM} = \frac{\rho_0}{5040}(\tilde{J}_1 M^{ij}_1 + \tilde{J}_2 M^{ij}_2 + \tilde{J}_3 M^{ij}_3 + \tilde{J}_4 M^{ij}_4) \ , \ (i,j=1,..,10) \tag{24}$$

Where $M^{ij}_k$ (k = 1, 2, 3, 4) are given by

$$M^{ij}_1 = \begin{pmatrix} 6 & 0 & 0 & 0 & 0 & -2 & 0 & 0 & -2 & -2 \\ 0 & 2 & 1 & 1 & -4 & -2 & -4 & -4 & -2 & -2 \\ 0 & 1 & 2 & 1 & -4 & -2 & -4 & -4 & -2 & -2 \\ 0 & 1 & 1 & 2 & -4 & -2 & -4 & -4 & -2 & -2 \\ 0 & -4 & -4 & -4 & 24 & 8 & 12 & 12 & 8 & 4 \\ -2 & -2 & -2 & -2 & 8 & 8 & 8 & 4 & 4 & 4 \\ 0 & -4 & -4 & -4 & 12 & 8 & 24 & 12 & 4 & 8 \\ 0 & -4 & -4 & -4 & 12 & 4 & 12 & 24 & 8 & 8 \\ -2 & -2 & -2 & -2 & 8 & 4 & 4 & 8 & 8 & 4 \\ -2 & -2 & -2 & -2 & 4 & 4 & 8 & 8 & 4 & 8 \end{pmatrix}, \quad M^{ij}_2 = \begin{pmatrix} 2 & 0 & 1 & 1 & -4 & -4 & -2 & -2 & -4 & -2 \\ 0 & 6 & 0 & 0 & 0 & 0 & -2 & -2 & 0 & -2 \\ 1 & 0 & 2 & 1 & -4 & -4 & -2 & -2 & -4 & -2 \\ 1 & 0 & 1 & 2 & -4 & -4 & -2 & -2 & -4 & -2 \\ -4 & 0 & -4 & -4 & 24 & 12 & 8 & 8 & 12 & 4 \\ -4 & 0 & -4 & -4 & 12 & 24 & 8 & 4 & 12 & 8 \\ -2 & -2 & -2 & -2 & 8 & 8 & 8 & 4 & 4 & 4 \\ -2 & -2 & -2 & -2 & 8 & 4 & 4 & 8 & 8 & 4 \\ -4 & 0 & -4 & -4 & 12 & 12 & 4 & 8 & 24 & 8 \\ -2 & -2 & -2 & -2 & 4 & 8 & 4 & 4 & 8 & 8 \end{pmatrix}$$

$$M^{ij}_3 = \begin{pmatrix} 2 & 1 & 0 & 1 & -2 & -4 & -4 & -2 & -2 & -4 \\ 1 & 2 & 0 & 1 & -2 & -4 & -4 & -2 & -2 & -4 \\ 0 & 0 & 6 & 0 & -2 & 0 & 0 & -2 & -2 & 0 \\ 1 & 1 & 0 & 2 & -2 & -4 & -4 & -2 & -2 & -4 \\ -2 & -2 & -2 & -2 & 8 & 8 & 8 & 4 & 4 & 4 \\ -4 & -4 & 0 & -4 & 8 & 24 & 12 & 4 & 8 & 12 \\ -4 & -4 & 0 & -4 & 8 & 12 & 24 & 8 & 4 & 12 \\ -2 & -2 & -2 & -2 & 4 & 4 & 8 & 8 & 4 & 8 \\ -2 & -2 & -2 & -2 & 4 & 8 & 4 & 4 & 8 & 8 \\ -4 & -4 & 0 & -4 & 4 & 12 & 12 & 8 & 8 & 24 \end{pmatrix}, \quad M^{ij}_4 = \begin{pmatrix} 2 & 1 & 1 & 0 & -2 & -2 & -2 & -4 & -4 & -4 \\ 1 & 2 & 1 & 0 & -2 & -2 & -2 & -4 & -4 & -4 \\ 1 & 1 & 2 & 0 & -2 & -2 & -2 & -4 & -4 & -4 \\ 0 & 0 & 0 & 6 & -2 & -2 & -2 & 0 & 0 & 0 \\ -2 & -2 & -2 & -2 & 8 & 4 & 4 & 8 & 8 & 4 \\ -2 & -2 & -2 & -2 & 4 & 8 & 4 & 4 & 8 & 8 \\ -2 & -2 & -2 & -2 & 4 & 4 & 8 & 8 & 4 & 8 \\ -4 & -4 & -4 & 0 & 8 & 4 & 8 & 24 & 12 & 12 \\ -4 & -4 & -4 & 0 & 8 & 8 & 4 & 12 & 24 & 12 \\ -4 & -4 & -4 & 0 & 4 & 8 & 8 & 12 & 12 & 24 \end{pmatrix}$$

(25)

Finally, we use quadratic approximation to model the metric in the element domain. Quadratic metric - QM is given by

$$J \approx \varphi^i \tilde{J}_i, \quad (i = 1,..,10)$$

$$\tilde{J}_i = \left| \tilde{J}^i_{mn} \right|, \quad \tilde{J}^i_{mn} = J_{mn} \big|_{\text{node } i}, \quad (m,n = 1,2,3) \tag{26}$$

Where $\tilde{J}_i(X_{kI})$ stand for metric evaluated at the node $(i=1,..,10)$, and $\tilde{J}^i_{mn}$ is the jacobian matrix (14) at the node $(i=1,..,10)$. The QM mass matrix follows from the above QM assumption (26) and analytic integration (6) of (7)

$$M^{ij}_{QM} = \frac{\rho_0}{45360} \tilde{J}_r M^{ij}_r, \quad (r = 1,..,10) \tag{27}$$

Where the first and tenth $M^{ij}_r$ $(r = 1,..,10)$ are given by

$$M^{ij}_1 = \begin{pmatrix} 18 & -6 & -6 & -6 & 24 & 12 & 24 & 24 & 12 & 12 \\ -6 & -6 & -1 & -1 & 0 & 6 & 6 & 6 & 6 & 8 \\ -6 & -1 & -6 & -1 & 6 & 6 & 0 & 6 & 8 & 6 \\ -6 & -1 & -1 & -6 & 6 & 8 & 6 & 0 & 6 & 6 \\ 24 & 0 & 6 & 6 & -24 & -24 & -12 & -12 & -24 & -12 \\ 12 & 6 & 6 & 8 & -24 & -40 & -24 & -12 & -20 & -20 \\ 24 & 6 & 0 & 6 & -12 & -24 & -24 & -12 & -12 & -24 \\ 24 & 6 & 6 & 0 & -12 & -12 & -12 & -24 & -24 & -24 \\ 12 & 6 & 8 & 6 & -24 & -20 & -12 & -24 & -40 & -20 \\ 12 & 8 & 6 & 6 & -12 & -20 & -24 & -24 & -20 & -40 \end{pmatrix}, ..., M^{ij}_{10} = \begin{pmatrix} 12 & 8 & 6 & 6 & -12 & -20 & -24 & -24 & -20 & -40 \\ 8 & 12 & 6 & 6 & -12 & -24 & -20 & -20 & -24 & -40 \\ 6 & 6 & 24 & 0 & -12 & -12 & -12 & -24 & -24 & -24 \\ 6 & 6 & 0 & 24 & -12 & -24 & -24 & -12 & -12 & -24 \\ -12 & -12 & -12 & -12 & 32 & 32 & 32 & 32 & 32 & 32 \\ -20 & -24 & -12 & -24 & 32 & 96 & 48 & 32 & 64 & 96 \\ -24 & -20 & -12 & -24 & 32 & 48 & 96 & 64 & 32 & 96 \\ -24 & -20 & -24 & -12 & 32 & 32 & 64 & 96 & 48 & 96 \\ -20 & -24 & -24 & -12 & 32 & 64 & 32 & 48 & 96 & 96 \\ -40 & -40 & -24 & -24 & 32 & 96 & 96 & 96 & 96 & 288 \end{pmatrix} \tag{28}$$

Details of jacobian matrices $\tilde{J}^i_{mn}$ $(i = 5,..,10)$ and details of constant matrices $M^{ij}_r$ $(r = 2,..,8)$ are omitted to keep the text short. Please contact the corresponding author if needed.

QM mass matrix require 10 metric evaluations, 10 additive terms exist in (27), whereas for numerical integration using fifteen Gauss points $n_p = 15$, fifteen metric evaluations are required. QM mass matrix (27) and numerical integration $n_p = 15$ (9) have similar form, yet, (9) involve additional terms, roughly 50% more computations. In the next section it will be shown that both schemes have equivalent accuracy.

## 5. Preliminary numerical study.

If the generated mesh is a straight-sided mesh then (12) is an exact consistent mass matrix. However, in the case of curved-sided elements, the accuracy of the developed CM, LM and QM mass matrices should be studied numerically and compared to the performance of widely accepted schemes.

We consider the next element family

$$X_{11}=0 \quad, X_{21}=0 \quad, X_{31}=0 \quad, X_{12}=1 \quad, X_{22}=0 \quad, X_{32}=0$$
$$X_{13}=0 \quad, X_{23}=1 \quad, X_{33}=0 \quad, X_{14}=0 \quad, X_{24}=0 \quad, X_{34}=1$$
$$X_{15}=\frac{1}{2}+\varepsilon, X_{25}=0+\varepsilon, X_{35}=0+\varepsilon, X_{16}=\frac{1}{2}+\varepsilon, X_{26}=\frac{1}{2}+\varepsilon, X_{36}=0+\varepsilon \quad (29)$$
$$X_{17}=0+\varepsilon, X_{27}=\frac{1}{2}+\varepsilon, X_{37}=0+\varepsilon, X_{18}=0+\varepsilon, X_{28}=0+\varepsilon, X_{38}=\frac{1}{2}+\varepsilon$$
$$X_{19}=\frac{1}{2}+\varepsilon, X_{29}=0+\varepsilon, X_{39}=\frac{1}{2}+\varepsilon, X_{110}=0+\varepsilon, X_{210}=\frac{1}{2}+\varepsilon, X_{310}=\frac{1}{2}+\varepsilon$$

For $\varepsilon=0$ the above result in a simple tetrahedron with flat faces, edges between node 1-2,1-3 and 1-4 have length 1 and perpendicular to each other, nodes 5-10 are located exactly in the middle of their edges. For $\varepsilon \neq 0$ (29) return varying metric elements where $\varepsilon$ is a displacement of the node in the associate direction. Herein $\varepsilon$ is a random number uniformly distributed in the range $-\delta \leq \varepsilon \leq +\delta$, different $\varepsilon$ is generated for every component. In that sense, $\delta$ is the coarseness of the mesh, for $\delta=0$ it follows that $\varepsilon=0$ namely elements are simple, whereas $\delta=0.2$ lead to a coarse curved-sided mesh.

The next values of $\delta$ have been considered $\delta \in (0.0, 0.025, 0.05, 0.075, 0.1, 0.0125, 0.15, 0.175)$. For every value of $\delta$, 100 different elements have been defined, and mass matrix components have been computed using various schemes, comparison to exact values have been performed and an averaged absolute errors are reported.

Figure 1, illustrates the averaged absolute errors of CM, LM and QM mass matrices as a function of the mesh coarseness $\delta$. According to our expectations, LM is more accurate than CM and QM is more accurate than LM. Averaged absolute error grows as the mesh become coarser.

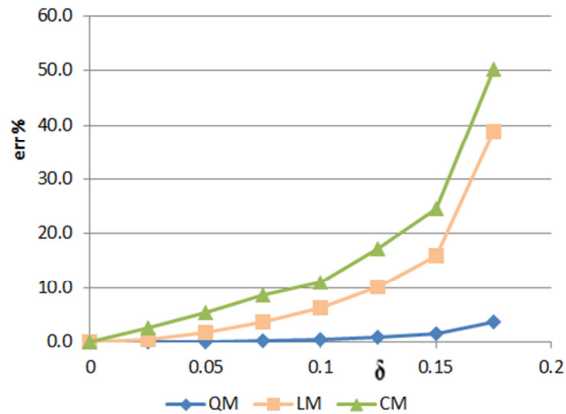

Figure 1: Averaged absolute error of CM, LM and QM semi-analytical mass matrices are reported as a function of the mesh coarseness $\delta$. For each point 100 random elements were used.

Figure 2, illustrates the averaged absolute error of the proposed CM rule vs numerically integrated mass matrix using four and five point quadrature. It is evident that our CM semi-analytical mass matrix significantly over performs in terms of accuracy numerical four and five

point schemes. As we stated earlier, mass matrix based on four point numerical integration, requires roughly four times more computations than CM, and five point integration is about five times more expensive. Moreover, as we stated earlier, CM is exact for straight sided element, figure 2 shows that an averaged absolute error for CM at $\delta = 0$ vanish.

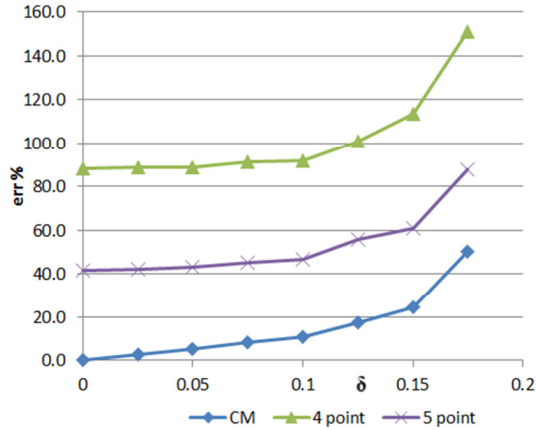

Figure 2: Averaged absolute error of CM semi-analytical rule and mass matrices based on four and five pint numerical integration schemes, presented as a function of the mesh coarseness $\delta$. For each point 100 random elements are used.

Figure 3, illustrates an averaged absolute error of QM and of mass matrix based on fifteen point numerical integration scheme. QM requires 10 evaluation of the metric (at nodal points) while fifteen point scheme requires 15 evaluations of the metric but rather equivalent accuracy is obtained.

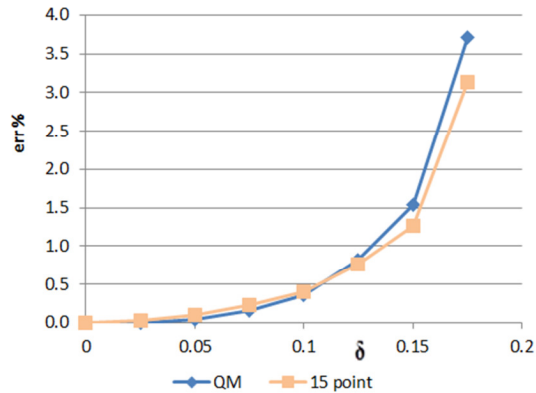

Figure 3: Averaged absolute error of QM semi-analytical mass matrix vs fifteen point numerical integration scheme as a function of the mesh coarseness $\delta$. For each point 100 random elements were used.

## 6. Conclusions

An exact analytical mass matrix for straight sided ten node tetrahedral element is derived and presented in a simple easy to implement form. Only one metric evaluation is needed, which makes the exact consistent mass matrix computationally inexpensive as the mass matrix obtained using one Gauss point numerical scheme, though considerably more accurate. Also, an exact analytical mass matrix for a curved-sided element is developed.

For a practical implementations of a curved-sided element, three systematic approximations are considered; constant metric (CM), linear metric (LM) and quadratic metric (QM). CM requires one metric evaluation at the centroid, LM requires four metric evaluations at the corner nodes 1-4 and QM requires ten metric evaluations performed at the nodes. Analytical integration is performed using the approximate metric models resulting in three closed-form mass matrices.

Preliminary numerical study using randomly generated coarse mesh is conducted. Our CM mass matrix is easy to implement and equivalent in computational effort to one point numerical integration nevertheless significantly over-performs in accuracy even four and five point numerical integration schemes. Widely used special fifteen point numerical integration scheme requires roughly 50% more computations than our QM mass matrix, while similar accuracy is established. Our CM, LM, and QM are sufficiently accurate easy to implement and considerably inexpensive approximations to consistent mass matrix of a curved-sided ten node tetrahedral element.

# References

[1] H. Borouchaki, P. George, and S. Lo, "Optimal Delaunay point insertion," *International journal for numerical methods in engineering,* vol. 39, 1996.

[2] S. Lo, "Optimization of tetrahedral meshes based on element shape measures," *Computers & structures,* vol. 63, pp. 951-961, 1997.

[3] C. Lee and S. Lo, "AUTOMATIC ADAPTIVE 3-D FINITE ELEMENT REFINEMENT USING DIFFERENT-ORDER TETRAHEDRAL ELEMENTS," *International journal for numerical methods in engineering,* vol. 40, pp. 2195-2226, 1997.

[4] P. Hammer, O. Marlowe, and A. Stroud, "Numerical integration over simplexes and cones," *Mathematical Tables and Other Aids to Computation,* pp. 130-137, 1956.

[5] G. Strang and G. J. Fix, *An analysis of the finite element method* vol. 212: Prentice-Hall Englewood Cliffs, NJ, 1973.

[6] P. Shiakolas, R. Nambiar, K. Lawrence, and W. Rogers, "Closed-form stiffness matrices for the linear strain and quadratic strain tetrahedron finite elements," *Computers & structures,* vol. 45, pp. 237-242, 1992.

[7] P. Shiakolas, K. Lawrence, and R. Nambiar, "Closed-form expressions for the linear and quadratic strain tetrahedral finite elements," *Computers & structures,* vol. 50, pp. 743-747, 1994.

[8] M. A. Moetakef, K. L. Lawrence, S. P. Joshi, and P. S. Shiakolas, "Closed-form expressions for higher order electroelastic tetrahedral elements," *AIAA Journal,* vol. 33, pp. 136-142, 1995.

[9] S. E. McCaslin, P. S. Shiakolas, B. H. Dennis, and K. L. Lawrence, "Closed-form stiffness matrices for higher order tetrahedral finite elements," *Advances in Engineering Software,* vol. 44, pp. 75-79, 2012.

[10] S. E. McCaslin, P. S. Shiakolas, B. H. Dennis, and K. L. Lawrence, "A New Approach to Obtaining Closed-Form Solutions for Higher Order Tetrahedral Finite Elements Using Modern Computer Algebra Systems," in *ASME 2011 International Design Engineering Technical Conferences and Computers and Information in Engineering Conference*, 2011, pp. 225-231.

[11] E. Hanukah and B. Goldshtein, "A structural theory for a 3D isotropic linear-elastic finite body," *arXiv preprint arXiv:1207.6767,* 2012.

[12] E. Hanukah, "Development of a higher order closed-form model for isotropic hyperelastic cylindrical body, including small vibration problem," *arXiv preprint arXiv:1312.0083,* 2013.

[13] E. Hanukah, "A new closed-form model for isotropic elastic sphere including new solutions for the free vibrations problem," *arXiv preprint arXiv:1311.0741,* 2013.

[14] E. Hanukah, "Higher order closed-form model for isotropic hyperelastic spherical shell (3D solid)," *arXiv preprint arXiv:1401.0204,* 2013.

[15] E. Hanukah, "Semi-analytical mass matrix for 8-node brick element," *arXiv preprint arXiv:1410.3195,* 2014.

[16] G. Dhondt, *The finite element method for three-dimensional thermomechanical applications*: John Wiley & Sons, 2004.

[17] P. Wriggers, *Nonlinear finite element methods*: Springer, 2008.

[18] M. Jabareen, E. Hanukah, and M. Rubin, "A ten node tetrahedral Cosserat Point Element (CPE) for nonlinear isotropic elastic materials," *Computational Mechanics,* vol. 52, pp. 257-285, 2013.